March 1, 2014
\input amstex
\magnification 1200
\documentstyle{amsppt}
\loadbold
\vsize 22truecm
\NoBlackBoxes
\def\qed{\hfill$\square$}
\def\ju{\vskip .8truecm plus .1truecm minus .1truecm}
\def\mju{\vskip .6truecm plus .1truecm minus .1truecm}
\def\sju{\vskip .4truecm plus .1truecm minus .1truecm}
\def\lra{\longrightarrow}

\def\qq{\lq\lq}
\def\lqq{\qq}
\def\pr{^\prime}
\def\prr{^{\prime\prime}}
\def\stm{\setminus}

\def\se{\subseteq}


\def\shimply{\Rightarrow} 
\def\lan{\langle}
\def\ra{\rangle}
\def\NN{\Bbb N}

\def\al{\alpha}

\def\ka{\kappa}
\def\be{\beta}
\def\ga{\gamma}

\def\ese#1{\buildrel{#1}\over\lra}
\def\sst#1#2{\{ {#1}:{#2} \}}
\def\hom#1#2{\hbox{\rm Hom}\,({#1},{#2})}
\def\homc#1#2#3{\hbox{\rm Hom}_{#1}({#2},{#3})}
\def\noi{\noindent}
\def\prf{\noi{\it Proof\quad}}
\def\rmod{R\bold{Mod}}
\def\Ker{\hbox{\rm Ker\,}}
\def\ima{\hbox{\rm Im\,}}
\def\ov{\overline}
\define\homkc{Hom$_*$-$\kappa$-commuting }
\define\kprod#1#2{\prod_{i\in{#1}}^\ka{#2}_i}
\document

\mju 
\centerline{\bf Hom$_*$ commuting with filtered products\footnote{The first version of this paper, dated February 11, 2007,  was communicated to C. U. Jensen}}

\centerline{\it Radoslav M. Dimitri\'c}
\vskip 2truecm

{\bf Abstract:} 
In a sufficiently rich category, such as a category of $R$-modules, 
and a given infinite cardinal $\ka$, we examine classes $\Cal H^\ka_*$ of objects $M$, such that a natural monomorphism is an isomorphism:

$${\prod_{i\in I}}^{\ka}\hom{M}{A_i}\cong\hom{M}{{\prod_{i\in I}}^{\ka} A_i},$$
for every family of objects $\sst{A_i}{i\in I}$ ($\prod^\ka$ denotes the subproduct of all vectors with support $<\ka$).

\ju
\S 1. {\bf Preliminaries. }

We will assume the axiom of choice or equivalently that every set may be well ordered; one consequence is the existence of arbitrary infinite products in the category of sets. Furthermore we will assume that the categories we work with have arbitrary products and coproducts. We will identify any cardinal $\ka$ with the smallest (initial) ordinal of that cardinality, when it is convenient to do so; this is equivalent to the statement that every cardinal is of the form $\aleph_\al$, for some ordinal number $\al$; we will use `fin' to denote a finite cardinal. Thus, an arbitrary non-empty index set $I$ may be assumed to be well ordered by an (initial) ordinal; if needed, we will assume ordinals to be regular, i.e. that cf\,$I$=$I$. A cardinal $\ka$ is regular, if it is not singular, i.e., if it cannot be represented as $\sup\sst{\al_i}{i<\theta}$, where each $\al_i,\ka$ and $\theta<\ka$ is a limit ordinal. Equivalently, $\ka$ is regular iff it cannot be represented as the sum of less than $\ka$ smaller cardinals. $|I|^+$ will denote the successor cardinal to cardinal $|I|$. An infinite cardinal $\aleph_\al$ is a successor cardinal, if it is of the form $\aleph_{\be+1}=\aleph_\al$; if $\al$ is a limit ordinal, then $\aleph_\al$ is called a limit cardinal.  Every successor ordinal is regular, but it is not always the case with limit ordinals, which may be singular. In this study, we are mostly interested in regular ordinals (cardinals).

If $a=(a_i)_{i\in I}\in\prod_{i\in I}A_i$, we will also write $a$ as a formal sum $\sum_{i\in I}a_i$ which is, more precisely, the formal sum $\sum_{i\in I}p_ia_i$, where $p_i:A_i\lra\prod A_i$ are the natural product injections. In absence of a topology/metric, the sum will \lqq make sense" when there are only totality of finitely many non-zero coordinates, at every coordinate of all $I$-vectors being summed. This will always be the case if we are doing genuine summation, not just the formal one. 
\vfill
\eject

In essence, our note concerns the category of unital (one-sided) $R$-modules, but we are using the language of general categories to indicate that the results and the proofs carry over to this more general setting, mutatis mutandis.

For an arbitrary family $\{ A_i, i\in I\}$ of (non-zero) objects, 
 and an arbitrary infinite cardinal $\ka$, ${\prod^\ka_{i\in I}} A_i$ will  denote the filtered $\ka$-product, namely it consists of all the vectors with support $<\ka$. 
The natural $\ka$-product-to-product ($\ka$ptp) embedding will be denoted by 
$u_\ka$ or $u:\kprod{I}{A}\lra\prod_{i\in I}A_i$.  

For $\ka\geq|I|^+$ ($\ka=\aleph_0$)  we have respectively 

$$
{\prod_{i\in I}}^\ka A_i=\prod_{i\in I}A_i \qquad \left({\prod_{i\in I}}^\ka A_i=\coprod_{i\in I}A_i\right). 
$$

{\bf Fact}: For every object $M$, we have a natural monomorphism 
$$\phi:{\prod_{i\in I}}^{\ka}\hom{M}{A_i}\lra\hom{M}{{\prod_{i\in I}}^{\ka} A_i}\eqno(*)$$
given by $(f_i:M\lra A_i)_{i\in I}\mapsto f=\sum_{i\in I}p_if_i:M\lra\prod_{i\in I}^\ka A_i$ (with coordinates $\forall\,\, \pi_if=f_i)$. This monomorphism is an isomorphism in case $\ka\geq |I|^+$, as well as when $I$ is a finite index set, regardless of $\ka$, thus, we can assume  in the sequel, when needed, without loss of generality,  that $\aleph_0\leq \ka\leq |I|$. 
 When $|I|$ is an infinite cardinal and $\ka\leq|I|$, then we do not necessarily have an isomorphism. For instance, if $\ka=\aleph_0$, (*) is not an isomorphism, for every $M$ and every family of non-zero modules $A_i$, even when all $A_i=R$. Thus we have the following 
\sju
{\bf Task 1.} Investigate, for various infinite cardinals $\ka$, and if possible characterize, objects $M$ 
such that for every infinite index set $I$, 
every family of (non-zero) objects $\{ A_i, i\in I\}$ monomorphism (*) is an isomorphism. Call every such object $M$ a {\it Hom$_*$-$\ka$-commuting object}. Given an infinite cardinal $\kappa$, denote by $\Cal H^\ka_*$ the class of all \homkc objects. For $\ka=\aleph_0$ the \homkc object was introduced by Mitchell, no later than 1965, under the name {\it small object} (see e.g also Rentchler, 1969 where small objects go under the name of $\Sigma$-{\it type} object (or $\Sigma$-generated) object). 
\sju

\S 2. {\bf Equivalent definitions.} 

While this is an ambitious task, we show in the sequel how to arrive at a number of illuminating results. We begin with the following:

\proclaim{Lemma 1}
Given an infinite cardinal $\ka$, an additive category $\Cal C$ with arbitrary coproducts and products, 
an arbitrary non-empty index set $I$  
  and a morphism $f:M\lra \kprod{I}{A}$, the following statements are equivalent:
\roster
\item The morphism $f$ factors through a $\ka\pr$- subproduct, $\ka\pr<\ka$, namely
there is a subset $J\se I$, $|J|=\ka\pr<\ka$ and a factorization $M\ese{f\pr}\kprod{J}{A}\ese{p_J}\kprod{I}{A}$,
$f=p_Jf\pr$ (here $p_J$ is the natural embedding of smaller $\kappa$-product into the larger; note that $\kprod{J}{A}=\prod_{i\in J}A_i$, since $|J|<\ka$).
\item $f=\sum_{i\in J}p_if_i$, for some $J\se I$, $|J|<\ka$ and for some morphisms $f_i:M\lra A_i$, $i\in J$.
\item If, for some $J\se I, |J|<\ka$, $q_{I\stm J}:\kprod{I}{A}\lra\kprod{I\stm J}{A}$ is the canonical quotient map $=\pi_{I\stm J}|\kprod{I}{A}$, then   $q_{I\stm J}f=0$. 
\item There is a subset $J\se I$, $|J|<\ka$ such that $\pi_{I\stm J}uf=0$.
\endroster
\endproclaim

\prf
As before, denote by $\pi\pr_i$ and $p\pr_i$ the natural product projections and injections associated
with the product indexed by $J$, and likewise by $u\pr$ the corresponding $\ka$ptp-morphism. 
Note that 

\noi (*) $\forall i\in J\, p_Jp_i\pr=p_i$ and (**) $\forall i\in J\,\,\, \pi_i\pr u\pr=\pi_i up_J$.

The proof is as follows:

 (1)$\shimply$ (2): $f=p_Jf\pr$$=p_J(\sum_{i\in J}p\pr_i\pi\pr_iu\pr)f\pr=$ (by (*))
$=\sum_{i\in J}p_i\pi\pr_iu\pr f\pr=$ (by (**))
\newline 
$=\sum_{i\in J}p_i\pi_iup_Jf\pr $$=\sum_{i\in J}p_i\pi_iuf$;
denote $f_i=\pi_iuf$, for all $i\in J$.

(2)$\shimply$(3): $q_{I\stm J}f=q_{I\stm J}(\sum_{i\in J}p_if_i)$$= \sum_{i\in J}q_{I\stm J}p_if_i=0$. 

(3)$\shimply$(4):  If $u\prr$ denotes the $\ka$ptp map associated with the $\ka$product on the index set $I\stm J$, then the proof is established by noting that $u\prr q_{I\stm J}=\pi_{I\stm J}u$ and thus $\pi_{I\stm J}uf=u\prr q_{I\stm J}f=0$. In fact the same observation proves the reverse implication. 

(4)$\shimply$(1): Equality $\pi_{I\stm J}uf=0$ ensures $\ima f\se\Ker(\pi_{I\stm J}u)\se \kprod{J}{A}$ and this in turn ensures validity of (1) \qed

We extend this result as follows, by not assuming a priori that index sets are the same or that the components $A_i$ are the same, in each of the equivalent statements:

\proclaim{Proposition 2}
Given an additive category  $\Cal C$ with arbitrary ($\ka$-products and) coproducts and an object $M$, the following are equivalent:
\roster
\item For every non-empty index set $I$ and an arbitrary family of objects $\sst{A_i}{i\in I}$ in 
$\Cal C$,  every morphism $f:M\lra\kprod{I}{A}$ in $\Cal C$ factors through a $\ka\pr$- subproduct, $\ka\pr<\ka$, namely there is a subset $J\se I$, $|J|=\ka\pr<\ka$ and a factorization $M\ese{f\pr}\kprod{J}{A}\ese{p_J}\kprod{I}{A}$, with $f=p_Jf\pr$.
\item  For every non-empty index set $I$ and an arbitrary family of objects $\sst{A_i}{i\in I}$ in 
$\Cal C$ and every morphism $f:M\lra\kprod{I}{A}$ in $\Cal C$,  there is a $J\se I$, $|J|<\ka$, such that
$f=\sum_{i\in J}p_if_i$, for some morphisms $f_i:M\lra A_i$, $i\in J$.
\item For every non-empty index set $I$ and an arbitrary family of objects $\sst{A_i}{i\in I}$ in 
$\Cal C$ and every morphism $f:M\lra\kprod{I}{A}$ in $\Cal C$ there is a $J\se I$, $|J|<\ka$, such that
$q_{I\stm J}f=0$, i.e. $\ima f\se\kprod{J}{A}$.
\item The functor $\homc{\Cal C}{M}{-}$ commutes with $\ka$-products, i.e., for every non-empty index set $I$
and an arbitrary family of objects $\sst{A_i}{i\in I}$ 
$$\phi:{\prod_{i\in I}}^{\ka}\hom{M}{A_i}\lra\hom{M}{{\prod_{i\in I}}^{\ka} A_i}\eqno(*)$$
via the natural isomorphism of Abelian groups 
$\phi:(f_i)_{i\in I} \mapsto \sum_{i\in I} p_if_i.$ 

\endroster
\endproclaim

\prf Lemma 1 establishes equivalence of the first three statements, since, a posteriori, it turns out that the index set $I$ and the product components $A_i$ may be the same in each of the equivalent statements. 

\noi (2)$\shimply$(4): Given a morphism $f:M\lra \kprod{I}{A}$ in $\Cal C$, there is a $J\se I$, $|J|<\ka$ such that  $f=\sum_{i\in J}p_if_i$, for some morphisms $f_i:M\lra A_i$. Note now that
$h=(f_{i0})_{i\in I}\in\prod_{i\in I}^\ka{\homc{\Cal C}{M}{A_i}}$ with $f_{i0}=f_i$, for $i\in J$ and $f_{i0}=0$ for $i\in I\stm J$ is such that $\phi(h)=f$.

\noi (4)$\shimply$(2): Given a morphism $f:M\lra \kprod{I}{A}$ in $\Cal C$, there is an 
$h=(f_i)_{i\in I}\in \prod_{i\in I}^\ka\homc{\Cal C}{M}{A_i}$ where, for some $J\se I$, $|J|<\ka$, $f_i=0$ whenever 
$i\in I\stm J$ and $\phi(h)=f$, i.e. $\sum_{i\in J}p_if_i=f$. 
 \qed

Finally, we have the following series of equivalent properties that could be used to define \homkc objects:

\proclaim{Theorem 3} In an additive category with infinite products and coproducts, given an infinite (regular) cardinal $\ka$, the following are equivalent, for an object $M$: 
\roster
\item for every infinite set $I$,
and every family of objects $\sst{A_i}{i\in I}$, the natural
monomorphism  defined in (*) is an isomorphism: 
$${\prod_{i\in I}}^{\ka}\hom{M}{A_i}\cong\hom{M}{{\prod_{i\in I}}^{\ka} A_i}\eqno(*I)$$
\newline(arbitrary Hom definition); 
\item $M$ is a \homkc with families of cardinality $\ka$, i.e., for every index set $J$ with $|J|=\ka$, 
and, for an arbitrary family of objects $\sst{A_i}{i\in J}$,
$${\prod_{i\in J}}^{\ka}\hom{M}{A_i}\cong\hom{M}{{\prod_{i\in J}}^{\ka} A_i}\eqno(*\ka)$$
via the natural isomorphism of Abelian groups 
\newline ($\ka$-Hom definition);

\item for an arbitrary non-empty index set $I$, for every family of objects
$\sst{A_i}{i\in I}$, and every morphism $f:M\lra \kprod{I}{A}$, there is a $J\pr\se I$ with $|J\pr|<\ka$, such that
for all $i\in I\stm J\pr$, $\pi_iuf=0$ 
\newline(arbitrary coordinatewise definition);

\item for every family $\{A_i, i\in J\}$, $|J|=\ka$,  of objects, and every morphism 
\newline $f:M\lra\kprod{J}{A}$, $\pi_iu_Jf=0$, for all 
$i\in J\stm J\pr$ for some $J\pr\se J$, $|J\pr|<\ka$. 
\newline ($\ka$-coordinatewise definition);

\item for every infinite index set $I$, 
for every family $\sst{A_i}{i\in I}$ of objects, and every
morphism $f:M\lra \kprod{I}{A}$, there is a well-ordering on $I$ such that, there is an $i_0\in I, i_0<\ka$
with $\pi_{i>i_0}f=0$ (arbitrary tailwise definition);

\item for every family $\{A_i\}_{i\in J}$ of objects, with $|J|=\ka$ ($J$ well-ordered), and every
morphism $f:M\lra \kprod{J}{A}$, there is an $i_0\in J$
with $i_0<\ka$ and $\pi_{i>i_0}f=0$ 
\newline ($\ka$-tailwise definition).
  \endroster
\endproclaim

(5) and (6) need regularity of $\ka$.

\prf Note that the index sets in (1), (3) and (5) are arbitrary, unlike (2), (4), (6) where index sets are of cardinality $\kappa$. Thus, each of (1),(3),(5) implies respectively (2),(4),(6). Equivalence of (1), (3) and (5) (arbitrary index sets) and equivalence of (2), (4) and (6) follow from Lemma 1 and Proposition 2. We only need to prove that one of the even numbered statements implies any of the odd numbered ones, to complete the proof of equivalence of all the statements. First, we show that regularity of $\ka$ is needed when working with tailwise definitions. 

(3)$\shimply$ (5): 
Assume that $f:M\lra \kprod{I}{A}$ is an arbitrary morphism. By (3), there is a $J\pr\se I$ with $|J\pr|<\ka$, such that for all $i\in I\stm J\pr$, $\pi_iuf=0$. We can well order $I$ in such a way that $|J\pr|$ is its initial segment (of cardinality $<\ka$).  Because, $\ka$ is regular, there is an $i_0<\ka$ in $I$ with $\pi_iuf=0$, for all $i>i_0$, which is the same as $\pi_{i>i_0}f=0$. This completes the proof of all the equivalences.

(4)$\shimply$ (6): The proof is the same, mutatis mutandis, as for 
(3)$\shimply$ (5).

(2)$\shimply$(4): Let $f:M\lra\kprod{J}{A}$ be an arbitrary morphism in $\Cal C$. By the assumption, isomorphism (*$\kappa$) holds, hence we can find morphisms $f_i:M\lra A_i$, with supp$(f_i)_{i\in J}<\ka$ i.e. there is a $J\pr\se J$ with $|J\pr|<\ka$ with $f_i=0$, for all $i\in J\stm J\pr$ and such that $f=\sum_{i\in J}p_if_i$ (summation of $<\ka$ non-zero summands, indexed by $J\pr$). Then $\forall\, i\in J\stm J\pr$, $\pi_iu_Jf=$ $ \pi_i\sum_{i\in J}u_Jp_if_i=$ $\pi_i(\sum_{i\in J\stm J\pr}up_if_i+\sum_{i\in J\pr}up_if_i)=0$.
 
(4)$\shimply$ (3): It is only non-trivial to consider cases when $|I|>\ka$. If, on the contrary, there is a $J\se I$ with $|J|=\ka$, such that $\forall i\in J$, $\pi_iuf\neq 0$, then we consider the cut of $f$:  $g=\pi_Jf:M\lra\kprod{J}{A}$. By (4), there exists a $J\pr\se J$ with $|J\pr|<\ka$ and $\forall i\in J\stm J\pr$, $\pi_iug\neq 0$, which would be a contradiction.     \qed
\sju

\S {\bf 3.  Examples, a characterization and constructions.}

 $gen\, M$ will denote the cardinality of a minimal set of generators of $M$ (and sometimes such a set of generators itself).\footnote{This practical notation was first introduced in Dimitric, 1984 and it has since been adopted in a number of cases, by other authors, without reference} 
and the contravariant Hom functors. 

\proclaim{Proposition 4} Let $\ka\leq |I|$ be an infinite cardinal
and let $M$ be an $R$-module  with $gen\, M<\ka$. Then the natural monomorphism
$$\phi:{\prod_{i\in I}}^{\ka}\hom{M}{A_i}\lra\hom{M}{{\prod_{i\in I}}^{\ka} A_i}\eqno(*)$$
is an isomorphism, for every family of objects $A_i, i\in I$, nemely $M$ is a \homkc object. Thus, if $gen\, M<\ka$, then $M\in \Cal H^{\ka\pr}_*$, for every $\ka\pr\geq\ka$.  In particular, every finitely generated $R$-module $M$ is Hom$_*-\ka$-commuting, for every infinite $\ka$. 
\endproclaim

\prf The reason for surjectivity is that $\forall f:M\lra\prod_{i\in I}^\ka A_i$, the image $gen\, \ima f\leq gen\, M<\ka$, and we may assume $gen\, f(M)=\{a^j=(a^j_i)_{i\in I}, j\in J\}$ with $|J|<\ka$ and support of every $a^j<\ka$. For every $i\in I$, we define $f_i:M\lra A_i$ as follows: For $m\in M$, let $f(m)=\sum_{j\in J_m}r_ja^j=$  (each $|J_m|\leq |J|$ is finite)
$\sum_{j\in J_m}r_j(a^j_i)_{i\in I}=\left(\sum_{j\in J_m}r_ja^j_i\right)_{i\in I}$ (finite sums); define $\forall i\in I$ $f_i(m)=\sum_{j\in J_m}r_ja^j_i$. Supp $(f_i)_{i\in I}<\ka$ since $\forall m\in M$ supp $f(m)<\ka$ and $|J|<\ka$. Clearly, by definition, $\phi(f_i)_{i\in I}=f$, which proves surjectivity. \qed

When $\ka=\aleph_0$, we get the fact that, if $M$ is a finitely generated object, then $\hom{M}{-}$ commutes with countable coproducts. 
\sju
{\bf Task 2.} Find and characterize categories (rings $R$) such that the only \homkc modules are those that are  $<\ka$-generated, as well as those where there are \homkc modules generated by at least $\ka$  elements. 

\proclaim{Proposition 5}
Let $\ka$ be any infinite limit cardinal and $M\in\rmod$ be such that 
{\bf no} ascending (smooth) $\ka$-chain of proper submodules of $M$ like 
$M_0<M_1<\dots< M_\al<\dots <M$, $\al<\rho\leq\ka$,  fills the whole of $M$, i.e. $\cup_{\al<\rho} M_\al=\sum M_\al\neq M$. Then $M$ is \homkc. 
\endproclaim

\prf By Theorem 3, it is sufficient to prove this for index sets $I$ with $\ka=|I|$; we well order $I$ so that $I$ represents the smallest ordinal of cardinality $\ka$. For every $\al<\ka$ we will denote by $\Pi_\al=\prod^{\al\ka}_{i\in I}A_i$ the truncated $\ka$-product that consists of elements $(a_i)_{i\in I}\in\Pi^\ka=\prod^{\ka}_{i\in I}A_i$ with $a_i=0$, for all $i>\al$.  Note that \{$\Pi_\al\}_{\al<\ka}$ is a smooth $\ka$-chain with union $\Pi^\ka$. 
Let $f\in \hom{M}{{\prod}^\ka_{i<\ka} A_i}$; denote by $M_\al=f^{-1}(\Pi_\al)$ 
(the inverse image of the $\al$-truncated $k$-product of $A_i$'s). Since $\{\Pi_\al\}_{\al<\ka}$ is a smooth $\ka$-chain uniting in $\Pi^\ka$, $\{M_\al\}_{\al<\ka}$ is likewise a smooth $\ka$-chain uniting in $M=f^{-1}(\Pi^\ka)$. By the assumption, this may happen only if not all the links are proper subobjects of $M$, i.e., if there exists an $\al<\ka$ with $M_\al=f^{-1}(\Pi_\al)=M$ This means that $f:M\lra\Pi_\al=\prod_{i\in I}^{\al\ka}A_i$. But  $\prod_{i\in I}^{\al\ka}A_i\cong\prod_{i\leq\al}A_i$ and, in this case $f=\Sigma_{i<\al}p_if_i$, for some $f_i:M\lra A_i$, which proves the statement. Note that the same proof goes if it is done by transfinite induction on $\ka$. \qed

\proclaim{Proposition 6} For an arbitrary infinite cardinal $\ka$, let $M$ be Hom$_*$-$\ka$-commuting. Then: 
No strictly ascending (smooth) $\rho$-chain $\rho\leq \ka$ of proper submodules of $M$ like 
$M_0<M_1<\dots< M_\al<\dots <M$, $\al<\rho\leq\ka$ fills the whole of $M$, i.e. $\cup_{\al<\rho} M_\al=\sum M_\al\neq M$.

\endproclaim
\prf 
 (here $\rho=\ka$) Assume that $M$ is Hom$_*-\ka$-commuting  and suppose, that on the contrary, for some such chain, we have $\cup_{i<\ka} M_i=\sum M_i=M$. Define 
$f(x)=(x+M_i)_{i<\ka}$ and hope that $f$ is a morphism $f:M\lra\kprod{\ka}{M/M}$. For every $x\in M$, there is the smallest $i_x<\ka$ with $x\in M_{i_x}$ and then, $\forall i\pr\geq i_x$,\, $x\in M_{i\pr}$ which implies that $|$supp$f(x)|\leq |i_x|<\ka$, thus, indeed $f\in\hom{M}{\kprod{\ka}{M/M}}$. Then, by the assumption, $f=(f_i)_{i<\ka}$, for some $(f_i)_{i<\ka}\in\prod_{i<\ka}^\ka\hom{M}{M/M_i}$. This means that there exists a $\ka_1<\ka$, such that $\forall i>\ka_1$\,\, $f_i=0$. This means that $\forall i>\ka_1$\,\, $\forall x\in M$, $x\in M_i$. This is a contradiction, since it is assumed that all $M_i$ are proper submodules of $M$. Hence, our claim holds.  \qed

\sju

\proclaim{Corollary 7} Let $\ka$ be an infinite cardinal and $M\in\rmod$ such that gen\,$M=\ka$. Then,
\roster
\item  If $\ka$ is a limit cardinal, $M$ is not \homkc.
\item Furthermore, if $\ka$ is a limit cardinal, $gen\, M\leq\ka$ and $M$ is Hom$_*-\ka$-commuting, then $gen\, M<\ka$. 
\endroster
\endproclaim

\prf (1) The reason is as follows: As usual, $\ka$ denotes the initial ordinal representing cardinal $\ka$.  If gen\,$M=\sst{m_\al}{\al<\ka}$, denote $M_\al=\lan\sst{m_i}{i\leq\al}\ra$. Then $M_\al, \al<\ka$ is an ascending $\ka$-chain of proper subobjects of $M$, with $\cup_{\al<\ka}M_\al=M$, because $\ka$ is a limit cardinal. 
By Proposition 6, this would be impossible, if $M$ were to be \homkc.

(2) This is a reformulation of (1). \qed

In any Abelian category with products and coproducts, we have the following:

\proclaim{Proposition 8} Let $\ka$ be an infinite (regular) cardinal. Then:
zero object?
\roster
\item The $0$ object is \homkc , thus the class of \homkc  objects in any category with the zero object is non-empty. Moreover all the $<\ka$-generated objects are in $\Cal H^\ka_*$.
\item An object is in $\Cal H^\ka_*$ iff all its quotient objects are in $\Cal H^\ka_*$.
\item If $0\lra A\ese{\al}B\ese{\be}C\lra 0$  is an exact
sequence and $A, C$ are in $\Cal H^\ka_*$   objects, then $B$ is likewise in $\Cal H^\ka_*$.
\item For an (infinite) index set $I$, and an arbitrary family of objects $\sst{B_i}{i\in I}$, the $\ka$-product $M=\kprod{I}{B}$ is 
in $\Cal H^\ka_*$, iff only finitely many $B_i\neq 0$. 
\item Given a finite ascending sequence $0=A_0<A_1<\dots<A_n=C$ of 
subobjects of an object
$C$, such that all the factors $A_{i+1}/A_i$, $i=0,1,\dots,n-1$
are in $\Cal H^\ka_*$, then $C$ is likewise in $\Cal H^\ka_*$. The claim is no longer true, if the length of chain is infinite (such as of countable cofinality). 
\item Let $F:\Cal A_1\lra\Cal A_2$ be an equivalence of
two categories with $\ka$-products and $M\in\Cal A_1$. Then $M$ is \homkc\,\,  in
$\Cal A_1$ if and only if $F(M)$ is \homkc  in $\Cal A_2$. 
\endroster
\endproclaim 

\prf
(1) is trivial, in view of Proposition 4. 

(2) holds, because a morphism $f:M/D\lra\kprod{I}{A}$ gives rise to a morphism $fq:M\lra\kprod{I}{A}$ ($q$ is the canonical quotient map). If $M$ is \homkc  then, by Proposition 2,  $fq$ is expressible as a sum $\sum_{i\in J} p_if_i$, for some $|J|<\ka$ and some  morphisms $f_i:M\lra A_i$. Given $\ov x\in M/D$, define, $\forall i\in I,$ $\ov f_i:M/D\lra A_i$ by $\ov f_i(\ov x)=f_i(x)$. $\ov f_i$ are morphisms since $f_i$ and $q$ are; moreover every $\ov f_i$ is well-defined, for if $\ov x=\ov x\pr$, then $f(\ov x)=fq(x)=fq(x\pr)=f(\ov x\pr)$ hence $\sum p_if_i(x)=\sum p_if_i(x\pr)$ and thus $\sum p_if_i(x-x\pr)=0$, $i\in J$. The latter is in $\kprod{I}{A}$, since $|J|<\ka$ thus every $p_if_i(x-x\pr)=0$ and this is possible only if for every $i\in I$,  $f_i(x)=f_i(x\pr)$. We can now represent $f=\sum p_i\ov f_i$ as a $J$-sum which, by Proposition 2, means that the quotient $M/D$ is \homkc. The other implication is trivial (once we take $D=0$).

(3) (Draw the commutative diagram to follow the argument easier). Consider an arbitrary morphism $f:B\lra\kprod{I}{A}\in\Cal C$. Then $f\al:A\lra\kprod{I}{A}$ and $A$ is assumed to be in $\Cal H^\ka_*$, hence there is a $J\se I$, $|J|<|I|$ such that $q_{I\stm J}f\al=0$, by Proposition 2.(3). 
Denote $\chi=q_{I\stm J}f$. By the universal property of the quotient (cokernel) construction, there is a unique $\ga:C\lra\kprod{I\stm J}{A}$ such that $\ga\be=\chi$. However $C$ is assumed to be in $\Cal H^\ka_*$, hence, there is a  $J\pr\se I\stm J$, $|J\pr|<\ka$ with $q\pr_{I\stm J\stm J\pr}\ga=0$, where $q\pr_{I\stm J\stm J\pr}:\kprod{I\stm J}{A}\lra \kprod{I\stm J\stm J\pr}{A}$ is the map corresponding to $\ga$, via Proposition 2.(3).  This implies $0=q\pr_{I\stm J\stm J\pr}\ga\be=q\pr_{I\stm J\stm J\pr}q_{I\stm J}f=q_{I\stm(J\cup J\pr)}f$. Since both $|J|, |J\pr|<\ka$, so is their union and  $q_{I\stm (J\cup J\pr)} f=0$, which, by Proposition 2.(3) establishes the fact that $B$ is in $\Cal H^\ka_*$. 

(4) If only finitely many components are $\neq 0$, $M$ is clearly \homkc, by way of canonical isomorphisms for finite products/coproducts. Assume for a moment that there are infinitely many $B_i\neq 0$ (say countably many). Then $M_\NN=\kprod{\NN}{B}$ is the ascending union of its proper subobjects $M_n=\kprod{n}{B}$, $n\in\NN$. By Proposition 6 this means that $\kprod{\NN}{B}$ is not \homkc. But then $\kprod{I}{B}$ cannot be \homkc, for otherwise, its quotient $M_\NN$ would have to be such as well, but this is impossible by (2) of this proposition. 

(5) The first claim follows from statements (0-3) of this proposition; the second claim follows from the fact that (say for a countable index set) a countable direct sum of \homkc objects is the union of an ascending chain of its proper subobjects (and by Proposition 6, not \homkc), 

(6) This statement is fully categorical and is straightforward. 
\sju

Our present effort concentrates on 
\sju

{\bf Task 3.} Give fairly detailed account of how the classes $\Cal H^\ka_*$ and $\Cal H^{\ka\pr}_*$ relate, for different cardinals $\ka$ and $\ka\pr$.
\ju
\ju
{\bf References}
\sju
{Dimitric, Radoslav}:
On pure submodules of free modules and $\kappa$-free modules. {\it CISM Courses and Lectures}, No. {\bf 287}, {\it Abeli
an Groups and Modules}, 373--381, Springer-Verlag, New York, 1984. 

{Mitchell, Barry}:
 {\it Theory of Categories} Academic Press, New York, 1965

{Rentschler, Rudolf M.}:
{Sur les modules $M$ tels que $\hom{M}{-}$ commute avec les
sommes directes,} {\it  C.R. Acad. Sc. Paris, S\'erie
A, }{\bf 268}{(1969), }{}{930--933}



\enddocument